\documentstyle[12pt]{article}
\begin{document}
\title{ HOMOLOGICAL DIMENSION OF  CROSSED PRODUCTS
\thanks {This work is supported by National Science Foundation}}
\author{ Shouchuan Zhang
\\ Department  of Mathematics,
Nanjing    University, 210008 } 
\date{}
\newtheorem{Theorem}{\quad Theorem}[section]
\newtheorem{Proposition}[Theorem]{\quad Proposition}
\newtheorem{Definition}[Theorem]{\quad Definition}
\newtheorem{Corollary}[Theorem]{\quad Corollary}
\newtheorem{Lemma}[Theorem]{\quad Lemma}
\newtheorem{Example}[Theorem]{\quad Example}
\maketitle


\addtocounter{section}{-1}

\noindent

\section {Introduction}

Throughout this paper, $k$ is a field, $R$  is an algebra
over $k$, and  $H$ is a Hopf algebra over $k$ .
We say that
$R \#_\sigma H$  is the crossed product of $R$ and $H$ if
$R\#_\sigma H$ becomes an algebra over $k$ by multiplication:
$$(a \# h) (b \# g) =
\sum _{h, g } a (h_1 \cdot b) \sigma (h_2, g_1) \# h_3g_2 $$
for any $a, b \in R, h, g \in H,$
where $\Delta (h) = \sum h_1\otimes h_2$
\ \ \ (  see, \cite [Definition 7.1.1] {Mo92}.)

Let $lpd ({}_RM)$, $lid({}_RM)$ \ \  and \ \  $lfd ({}_RM)$ denote the
left projective dimension, left injective dimension and left flat
dimension of left $R$-module $M$, respectively.
Let $lgD (R)$ \ \  and \ \  $wD (R)$ denote the
left global  dimension and weak dimension of algebra $R$, respectively.

Crossed products are very  important algebraic structures. The relation
between homological dimensions of algebra $R$  and crossed product
$R \#_\sigma H$ is often  studied. J.C.Mconnell and J.C.Robson in \cite
[Theorem 7.5.6] {MR} obtained that
$$rgD (R) = rgD (R*G)$$ for
any finite group $G$ with $\mid G \mid ^{-1} \in k$, where $R*G$ is skew group
ring. It is clear that every skew group ring
$R*G$  is a crossed product $R \# _\sigma kG$ with
trivial $\sigma$.
Zhong Yi
in \cite {Y}
 obtained            that
the global dimension of crossed product $R*G $ is finite when
the global dimension
of $R$ is finite  and some other conditions hold.


In this paper,
we obtain that the global dimensions of $R$ and
 the crossed product
$R \# _\sigma H$ are the same; meantime,
 their weak dimensions are also the same,  when 
 $H$ is finite-dimensional semisimple and cosemisimple Hopf algebra.

                \section {The homological dimensions of modules over 
                crossed products}

In this section, we give the relation between homological dimensions of
modules over $R$  and  $R\#_\sigma H$.
  
If $M$ is a left (right) $R\#_\sigma H$-module, then $M$ is also a left
(right ) $R$-module since we can view $R$ as a subalgebra of $R\#_\sigma H$.

\begin {Lemma} \label {1.1}
Let $R$  be a subalgebra of algebra $A$.

(i) If $M$ is a free $A$-module and $A$ is a free $R$-module, then
$M$  is   a free $R$-module;

(ii) If $P$ is a projective left $R\#_\sigma H$-module,
then $P$ is a projective left $R$-module;

(iii) If $P$  is a projective right $R\#_\sigma H$-module and
$H$ is a Hopf algebra with invertible antipode, then
$P$ is a projective right $R$ -module;

(iv) If $${\cal P}_M \hbox { \ \ \ \ } \hbox { \ \ \ \ }: \cdots  P_n \stackrel {d_n} {\rightarrow} P_{n-1} \cdots
\rightarrow P_0 \stackrel {d_0}{\rightarrow}  M \rightarrow 0 $$
 is  a projective resolution of left $R\#_\sigma H$-module $M$,
 then ${\cal P}_M$ is a projective resolution of left $R$-module $M$;

 (v) If
 $${\cal P}_M: \hbox { \ \ \ \ } \cdots  P_n \stackrel {d_n} {\rightarrow} P_{n-1} \cdots
\rightarrow P_0 \stackrel {d_0}{\rightarrow}  M \rightarrow 0 $$
 is a projective resolution of right  $R\#_\sigma H$-module $M$ and
 $H$ is a Hopf algebra with invertible antipode, then
  ${\cal P}_M$  is a projective resolution of right $R$-module $M$ .

\end {Lemma}

{\bf Proof.}  (i) It   is  obvious.

(ii)  Since $P$ is a projective $R\#_\sigma H$-module, we have that
there exists a free $R \#_\sigma H$-module $F$ such that
$P$ is a summand of $F$.
 It is clear that $R\# _\sigma H \cong R \otimes H$  as left $R$-module,
 which implies that  $R\#_\sigma H$ is a free $R$-module.
Thus it follows from part (i)  that $F$ is a free $R$-module and
$P$ is a summand of $F$ as $R$-module.
Consequently, $P$ is a projective  $R$-module.

(iii)  By \cite [Corollary 7.2.11] {Mo92}, $R\#_\sigma H \cong H \otimes R$
as right $R$-module. Thus $R \#_\sigma H$  is a free right $R$-module.
Using the  method  in the proof of part (i), we have that $P$ is a projective rigth $R$-module.

(iv) and (v) can be obtained by part (ii) and (iii).           
$\Box$

\begin {Lemma} \label {1.2.1}
(i) Let $R$ be a subalgebra of $A$.
If $M$ is a  flat right ( left ) $A$-module and $A$ is a flat right ( left ) $R$-module, then
$M$  is   a flat right ( left ) $R$-module;

(ii) If $F$ is a flat left $R\#_\sigma H$-module,
then $F$ is a flat left $R$-module;

(iii) If $F$  is a flat right $R\#_\sigma H$-module and
$H$ is a Hopf algebra with invertible antipode, then
$F$ is a flat right $R$ -module;

(iv) If $${\cal F}_M: \hbox { \ \ \ \ } \cdots  F_n \stackrel {d_n} {\rightarrow} F_{n-1} \cdots
\rightarrow F_0 \stackrel {d_0}{\rightarrow}  M \rightarrow 0 $$
 is  a flat resolution of left $R\#_\sigma H$-module $M$,
 then ${\cal F}_M$ is a flat resolution of left $R$-module $M$;

 (v) If
 $${\cal F}_M: \hbox { \ \ \ \ } \cdots  F_n \stackrel {d_n} {\rightarrow} F_{n-1} \cdots
\rightarrow F_0 \stackrel {d_0}{\rightarrow}  M \rightarrow 0 $$
 is a flat resolution of right  $R\#_\sigma H$-module $M$ and
 $H$ is a Hopf algebra with invertible antipode, then
  ${\cal F}_M$  is a flat resolution of $M$;

\end {Lemma}
{\bf Proof.} (i)  We only show part (i) in the case which 
$M$ is a right $A$-module and $A$ is a right $R$-module; the other cases
 can similarly    be shown.
Let $$0 \rightarrow X \stackrel {f}{\rightarrow} Y$$
be an exact left $R\#_\sigma H$-module sequence. By assuptions,
    $$0 \rightarrow A  \otimes _R X \stackrel { A \otimes f}{\rightarrow}
    A \otimes _R Y$$
    and 
    $$0 \rightarrow M \otimes _A (A \otimes _R X )\stackrel {M \otimes
    (A\otimes f)}
    {\rightarrow} M \otimes _A (A \otimes _R Y)$$
    are exact sequences.
    Obviously, $$M \otimes _A (A \otimes _R X) \cong M \otimes _RX
    \hbox { \ \  \ and }  M \otimes _A (A \otimes _R Y) \cong M \otimes _RY $$
    as additive groups.
    Thus 
    $$0 \rightarrow M  \otimes _R X \stackrel {M\otimes f}
    {\rightarrow} M  \otimes _R Y$$
    is an exact  sequence, which implies
    $M$ is a flat $R$-module.

    (ii)-(v)  are immediate consequence of part (i)
    $\Box$

              The following is a immediate consequence of Lemma \ref {1.1} and
             \ref {1.2.1}.
             
\begin {Proposition} \label {1.2}

(i) If $M$ is a left $R\#_\sigma H$-module,
then $$lpd ({}_RM) \leq lpd({}_{R \#_\sigma H}M);$$

(ii) If $M$ is a right $R\#_\sigma H$-module and
$H$ is a Hopf algebra with invertible antipode, then
 $$rpd (M_R) \leq rpd({M}_{R \#_\sigma H})$$

(iii) If $M$ is a left $R\#_\sigma H$-module,
then $$lfd ({}_RM) \leq lfd({}_{R \#_\sigma H}M);$$

(iv) If $M$ is a right $R\#_\sigma H$-module and
$H$ is a Hopf algebra with invertible antipode, then
 $$rfd (M_R) \leq rfd({M}_{R \#_\sigma H}).$$

\end {Proposition}

\begin {Lemma} \label {1.4}
 Let $H$ be a finite-dimensional semisimple  Hopf algebra,
 and let $M$ and $N$ be left $R\#_\sigma H$-modules.
  If
 $ f$  is an $R $-module homomorphism from $M$ to $N$,
  and
  $$\bar f (m) = \sum \gamma ^{-1}(t_1)f(\gamma (t_2)m)$$
 for any $m \in M,$ then 
 $\bar f$  is an $R \#_\sigma H$-module homomorphism from $M$ to $N$,
 where  $t \in \int _H^r$ with $\epsilon (t)=1$, and
 $\gamma $  is a map from $H$ to $R \#_\sigma H$ sending $h$ to $1 \# h.$

\end {Lemma}
{\bf Proof.} (see, the proof of \cite [Theorem 7.4.2] {Mo92})
 For any $a\in R, h \in H , m \in M$, we see that
\begin {eqnarray*}
\bar f (am) &=&\sum \gamma ^{-1}(t_1)f(\gamma (t_2)am ) \\
&=&\sum \gamma ^{-1}(t_1)f((t_2 \cdot a) \gamma (t_3)m ) \\
&=&\sum \gamma ^{-1}(t_1)(t_2 \cdot a)f(\gamma (t_3)m ) \\
 &=&\sum a \gamma ^{-1}(t_1)f(\gamma (t_2)m ) \\
 &=&a \bar f (m)
 \end {eqnarray*}
and
\begin {eqnarray*}
 \bar f (\gamma (h)m) &=&\sum \gamma ^{-1}(t_1)f(\gamma (t_2)
 \gamma (h)m ) \\
&=&\sum \gamma ^{-1}(t_1)f(\sigma (t_2, h_1) \gamma (t_3h_2)m )
\hbox {\ \ \ by \cite [Definition 7.1.1] {Mo92}} \\
&=&\sum \gamma ^{-1}(t_1)\sigma (t_2 , h_1)f(\gamma (t_3h_2)m ) \\
&=& \sum \gamma (h_1)\gamma ^{-1}(t_1h_2) f(\gamma (t_2h_3)m)\\
&=&\sum  \gamma (h) \gamma ^{-1}(t_1)f(\gamma (t_2)m )
 \hbox { \ \ \ since } \sum h_1 \otimes t_1h_2 \otimes t_2 h_3 = 
 \sum h \otimes t_1 \otimes t_2 \\
  &=&\gamma (h) \bar f (m)
 \end {eqnarray*}
   Thus $\bar f$ is an $R\#_\sigma H$-module homomorphism.
  $\Box$

  In fact, we can obtain a functor by Lemma \ref {1.4}. Let ${}_{R\#_\sigma 
  H} \overline {\cal M}$ denote the full subcategory of ${}_R {\cal M};$
  its objects are all of left $R\#_\sigma H$-modules and its morphisms 
from $M$  to $N$  are all of $R$-module homomorphisms from $M$ to $N$. 
For any  $M, N \in ob 
   {} _{R\#_\sigma H} \overline {\cal M}$  and $R$-module homomorphism 
   $f$ from $M$  to $N$, we define that 
   $$F:  {}_{R\#_\sigma H}\overline {\cal M} \longrightarrow 
   {}_ {R\#_\sigma H} {\cal M}$$
   such that $$F(M) = M \hbox   {  \ \ \ \ and \ \ \ \  } F(f) = \bar f,$$
   where $\bar f$ is defined in Lemma \ref {1.4}. It is clear that 
   $F$  is a functor.

\begin {Lemma} \label {1.5}
 Let $H$ be a finite-dimensional semisimple  Hopf algebra,
 and let $M$ and $N$ be right $R\#_\sigma H$-modules.
  If
 $ f$  is an $R$-module homomorphism from $M$ to $N$,
  and
  $$\bar f (m) = \sum f(m \gamma ^{-1}(t_1))\gamma (t_2)$$
\newpage 
 for any $m \in M,$ then 
 $\bar f$  is an $R \#_\sigma H$-module homomorphism from $M$ to $N$,
 where  $t \in \int _H^r$ with $\epsilon (t)=1,$
 $\gamma $  is a map from $H$ to $R \#_\sigma H$ sending $h$ to $1 \# h.$

\end {Lemma}
{\bf Proof.} (see, the proof of \cite [Theorem 7.4.2] {Mo92})
 For any $a\in R, h \in H , m \in M$, we see that
\begin {eqnarray*}
\bar f (ma) &=&\sum f(ma \gamma ^{-1}(t_1))\gamma (t_2)  \\
&=&\sum f(m \gamma ^{-1}(t_1)(t_2 \cdot a)) \gamma (t_3) \\
&=&\sum f(m \gamma ^{-1}(t_1))(t_2 \cdot a) \gamma (t_3) \\
 &=&\sum f(m \gamma ^{-1}(t_1))\gamma (t_2) a \\
 &=& \bar f (m) a
 \end {eqnarray*}
and
\begin {eqnarray*}
 \bar f (m \gamma (h)) &=&\sum f(m \gamma (h) \gamma ^{-1}(t_1))
 \gamma (t_2)  \\
&=&\sum f(m \gamma (h_1) \gamma ^{-1} (t_1 h_2)) \gamma (t_2h_3 )
 \hbox { \ \ \ since } \sum h_1 \otimes t_1h_2 \otimes t_2 h_3 =
 \sum h \otimes t_1 \otimes t_2 \\
&=&\sum f(m \gamma ^{-1}(t_1) \sigma (t_2, h_1)) \gamma (t_3h_2 )
\hbox {\ \ \ by \cite [Definition 7.1.1] {Mo92}} \\
&=&\sum f(m \gamma ^{-1}(t_1))\sigma (t_2 , h_1)\gamma (t_3h_2 )
  \\
 &=&\sum f(m \gamma ^{-1}(t_1)) \gamma (t_2)\gamma (h) \\
  &=& \bar f(m)) \gamma (h)
 \end {eqnarray*}
 Thus $\bar f$  is an $R\#_\sigma H$-module homomorphism.
  $\Box$

 \begin {Proposition} \label {1.6}  Let
 $H$ be a finite-dimensional semisimple Hopf algebra.

(i) If $P$  is a left (right) $R\#_\sigma H$-modules and a projective
left (right) $R$-module,
then $P$ is a projective left (right) $R \#_\sigma H$-module;

(ii) If $E$  is a left (right) $R\#_\sigma H$-modules and an injective
left (right) $R$-module,
then $E$ is an injective left (right) $R \#_\sigma H$-module;

(iii) If $F$  is a left (right) $R\#_\sigma H$-modules and a flat 
left (right) $R$-module,
then $F$ is a flat left (right) $R \#_\sigma H$-module.

\end {Proposition}
{\bf Proof.}
(i)  Let $$X \stackrel {f} {\rightarrow} Y \rightarrow 0$$
be an exact sequance of left (right ) $R\#_\sigma H$-modules and
$g$  be a $R\#_\sigma H$-module homomorphism from $P$ to $Y$. Since $P$
is a projective left (right) $R$-module, we have that
there exists a $R$-module homomorphism $\varphi $
from $P$ to $X$, such that  $$f \varphi = g.$$
By Lemma \ref {1.4} and \ref {1.5}, there exists a
$R\#_\sigma H$-module homomorphism  $\bar \varphi $   from
$P$  to $X$ such that  $$f \bar \varphi = g.$$ 
\newpage
Thus $P$ is a projective left (right )  $R\#_\sigma H$-module.

Similarly, we can obtain the proof of part (ii).

(iii) Since $F$  is a flat left (right ) $R$-module, we have the character module
$ Hom _{\cal Z}
(F, {\cal Q}/{\cal Z})$ of $F$ is a injective left (right )  $R$-module by \cite
[Theorem 2.3.6] {T}. Obviously, $Hom_{\cal Z} (F, {\cal Q}/{\cal Z})$ is
a left (right ) $R\#_\sigma H$-module. By part (ii),
 $Hom_{\cal Z} (F, {\cal Q}/{\cal Z})$ is a injective left ( right ) $R \#_\sigma H$
 -module.  Thus $F$ is a flat left (right ) $R\#_\sigma H$-module.
 $\Box$

\begin {Proposition} \label {1.7}
Let $H$ be a finite-dimensional semisimple Hopf algebra. Then
 for  left (right ) $R \#_\sigma H $ -modules $M$ and $N$,
 $$Ext_{R\#_\sigma H}^n (M, N) \subseteq Ext_R ^n (M,N), $$
 where $n$ is any natural number.

\end {Proposition}
{\bf Proof.} We view the $Ext^n(M,N)$ as the equivalent classes of $n$-
extension of $M$ and $N$ (see, \cite
[Definition 3.3.7] {T}). We only prove  this result for  $n=1$. 
For  other cases, we can
similarly prove.
We denote the equivalent classes in
$Ext_{R\#_\sigma H}^1 (M,N)$  and $Ext_R^1(M,N)$ by  [E] and
[F]' , respectively, where  $E$ is an extension of $R\#_\sigma H$-modules $M$ and $N$,
 and $F$ is an extension  of $R$-modules $M$ and $N$.
We define a map $$\Psi :Ext _{R\#_\sigma H}^1 (M,N) \rightarrow
Ext_R^1(M,N), \hbox { \ \ \ by sending  \ } [E] \hbox { \ \ to }  [E]' .$$
Obviously, $\Psi$ is a map. Now we show that $\Psi $  is injective.
Let
 $$0 \rightarrow M \stackrel {f}{\rightarrow}  E
 \stackrel {g}{\rightarrow } N \rightarrow 0 
\hbox {  \ \ \ and \ } 0 \rightarrow M \stackrel {f'} {\rightarrow} E'
  \stackrel {g'}{\rightarrow}  M \rightarrow 0 $$
  are two extensions of $R\#_\sigma H$-modules $M$ and $N$, and they
  are equivalent in $Ext_R^1 (M,N)$. Thus there exists
  $R$-module homomorphism $\varphi $ from $E$ to $E'$ such that
$$\varphi f =f' \hbox { \ \ \ and \ \ } \varphi g = g'.$$
  By lemma \ref {1.4} , there exists $R\#_\sigma H$-module homomorphism
  $\bar \varphi $  from $E$  to $E'$ such that  
   $$\bar \varphi f =f' \hbox { \ \ \ and \ \ } \bar \varphi g = g'.$$
Thus $E$ and $E'$  is equivalent in $Ext _{R\#_\sigma H}^1 (M,N),$
which implies that $\Psi$ is  injective. $\Box$

\begin {Lemma} \label {1.7}
 For any $M \in {\cal M}_{R \#_\sigma H} $  and $N \in {}_{R \#_\sigma H}
 {\cal M}$, there exists an additive group  homomorphism
 $$\xi:   M \otimes _{R} N  \rightarrow M \otimes _{R \#_\sigma H}N$$
    by sending  $(m \otimes n)$  to $m \otimes n$,
    where $m \in M, n \in N$.

\end {Lemma}
{\bf Proof.} It is trivial. $\Box$

 \begin {Proposition} \label {1.8}

 If $M$   is a right $R\#_\sigma H$-modules and
$N$ is a left    $R \#_\sigma H$-module,   
then
     there exists additive group homomorphism 
     $$\xi_* : {\ } Tor ^R_n (M,N)  \longrightarrow 
     Tor ^{R\#_\sigma H}_n (M,N)$$
     such that   $\xi_*([z_n] ) = [\xi(z_n)]$,
     where $\xi$ is the same as in Lemma \ref {1.7}.

\end {Proposition}
{\bf Proof.}
 Let
$${\cal P}_M: \hbox { \ \ \ \ } \cdots  P_n \stackrel {d_n} {\rightarrow} P_{n-1} \cdots
\rightarrow P_0 \stackrel {d_0}{\rightarrow}  M \rightarrow 0 $$
 is a projective resolution of right  $R\#_\sigma H$-module $M$,  and
 set
 $$T = -   {\otimes }_{R\#_\sigma H} N \hbox { \ \ \  and  \ \ }
 T^R = - {\otimes }_ {R} N .$$
We have that
  $$ T{\cal P}_{\hat M}: \hbox { \ \ \ \ } \cdots  T(P_n) \stackrel {Td_n} {\rightarrow}
  T(P_{n-1})
  \cdots
\rightarrow T(P_1) \stackrel {Td_1}{\rightarrow}  TP_0 \rightarrow 0 $$
and  
$$T^R{\cal P}_{ \hat  M}:\hbox  { \ \ \ } \cdots T^R( P_n) \stackrel {T^Rd_n} {\rightarrow}
T^R( P_{n-1}) \cdots
\rightarrow T^R(P_1) \stackrel {T^Rd_1}{\rightarrow} T^R (P_0) \rightarrow 0 $$
are complexes .
     Thus $\xi$ is a complex homomorphism from
     $T^R {\cal P}_{\hat M}$ to $T{\cal P}_  {\hat M}$,   
     which implies that
     $\xi_*$  is an additive group homomorphism.
     $\Box$

\section {The global dimensions and weak dimensions of crossed products}

In this section we give the relation between homological dimensions of $R$
 and  $R \#_\sigma H$.

\begin {Lemma} \label {2.2}  
If $R$ and $R'$  are Morita equivalent rings, then

(i) rgD(R) = rgD(R');

(ii) lgD(R) = lgD(R');

(iii) wD(R)= wD(R').
\end {Lemma}

{\bf Proof.}  It is an immediate consequence of \cite [Proposition 21.6,
Exercise  22.12] {AF} $\Box$

\begin {Theorem} \label {2.4} Let
If $H$  is a finite-dimensional semisimple Hopf algebra,

(i)  $rgD(R\#_\sigma H)\leq rgD(R) ;$

(ii)  $ lgD(R\#_\sigma H) \leq lgD(R);$

(iii) $ wD(R\# _\sigma H) \leq wD(R)$.

\end {Theorem}

{\bf Proof.} 
(i)
When
$lgD(R)$ is infinite, obviously  part (i) holds.
Now we assume $lgD(R) =n.$
For any left $R\#_\sigma H$-module $M,$ and
 a projective resolution of left  $R\#_\sigma H$-module $M$ :
  $${\cal P}_M: \hbox { \ \ \ \ } \hbox { \ \ \ \ } \cdots  P_n \stackrel {d_n} {\rightarrow} P_{n-1} \cdots
\rightarrow P_0 \stackrel {d_0}{\rightarrow}  M \rightarrow 0 ,$$
we have  that ${\cal P} _M$ is also  a projective resolution of
left  $R$-module $M$ by Lemma \ref {1.1}.
Let $K_n = ker { \ } d_n$  be syzygy $n$ of ${\cal P}_M$. Since $
lgD(R) = n$, $Ext _R^{n+1} (M, N)=0$ for any left $R$-module
$N$  by \cite [Corollary 3.3.6] {T}. Thus $Ext _R^1 (K_n , N)=0,$  which implies
$K_n$ is a projective $R$-module. By Lemma \ref {1.6} (i),
$K_n$ is a projective $R\#_\sigma H$-module and
$Ext_{R\#_\sigma H}^{n+1}(M,N)= 0$  for any $R\#_\sigma H$-module $N$.
Consequently, $$lgD(R\#_\sigma H) \leq n = lgD(R) \hbox { \ \ \ by
\cite [Corollary 3.3.6] {T} }.$$

We complete the proof of part (i).

We can similarly show part (ii)  and part (iii). $\Box$

\begin {Theorem} \label {2.5}  Let
$H$  be a finite-dimensional semisimple and cosemisimple Hopf algebra.
Then

(i) $rgD(R) = rgD(R\#_\sigma H);$

(ii) $rgD(R) = rgD(R\#_\sigma H);$

(iii) $wD(R)= wD(R\# _\sigma H)$.

\end {Theorem}

{\bf Proof.} 
(i)  By dual theorem (see, \cite [Corollary 9.4.17] {Mo92}), we have
$(R \#_\sigma H)\# H^*$  and $R$  are Morita equivanlent  algebras.
Thus $lgD(R) =  lgD ((R \#_\sigma H)\#H^*)$ by Lemma \ref {2.2} (i).
Considering Theorem \ref {2.4} (i), we have that
$$ lgD ((R \#_\sigma H)\#H^*) \leq lgD (R \#_\sigma H) \leq lgD (R) .$$
Consequently,
$$lgD (R) = lgD (R\#_\sigma H).$$

Similarly, we can prove (ii) and (iii) .
$\Box$

\begin {Corollary} \label {2.6}  Let
$H$  be a finite-dimensional semisimple  Hopf algebra.

(i) If $R $ left (right ) semi-hereditary, then  so is $R\#_\sigma H;$

(ii) If $R $  is von Neumann regular,  then  so is $R\#_\sigma H.$

\end {Corollary}

{\bf Proof.} (i)  It follows from Theorem \ref {2.4} and
\cite [Theorem 2.2.9] {T}          .

(ii) It follows from Theorem \ref {2.4}  and \cite [Theorem 3.4.13] {T}.
     $\Box$

By the way, part (ii) of Corollary \ref {2.6} give one case about the 
semiprime question in \cite [Question 7.4.9] {Mo92}. That is, 
If $H$ is a finite-dimensional  semisimple Hopf algebra and  
$R$ is a von Neumann regular  algebra (notice that every von Neumann
regular algebra is  semiprime ),  then    $R\#_\sigma H$ is semiprime.

\begin {Corollary} \label {2.8}  Let
$H$  be a finite-dimensional semisimple  and cosemisimple Hopf algebra.
Then
 
 (i) $R$  is  semisimple artinian iff $R\#_\sigma H$  \ \ is semisimple artinian;

(ii) $R $ is left (right ) semi-hereditary iff  $R\#_\sigma H$ \ \  is
left (right ) semi-hereditary;

(iii)  $R $  is von Neumann regular iff $R\#_\sigma H $ \ \ is von Neumann regular.

\end {Corollary}
 
{\bf Proof.} 
(i) It follows from Theorem \ref {2.5}  and \cite [Theorem 2.2.9] {T}.

(ii)  It follows from Theorem \ref {2.5} and
\cite [Theorem 2.2.9] {T}.

(iii) It follows from Theorem \ref {2.5}  and \cite [Theorem 3.4.13] {T}.
$\Box$
           
If $H$ is commutative or cocommutative, then $S^2 = id_H$  by \cite {Sw}.
Consequently, by \cite [Proposition 2 (c)] {R94},
 $H$  is semisimple and cosemisimple iff
the character $char k $ of $k$ does  not
divides $dim H$. Considering Theorem \ref {2.5} and Corollary \ref {2.8},
we have:

\begin {Corollary} \label {2.9}

Let $H$ be a finite-dimensional commutative or cocommutative Hopf algebra.
If the character $char k $ of $k$ does  not
divides $dim H$, then

(i) $rgD(R) = rgD(R\#_\sigma H);$

(ii) $rgD(R) = rgD(R\#_\sigma H);$

(iii) $wD(R)= wD(R\# _\sigma H);$
 
 (iv) $R$  is  semisimple artinian iff $R\#_\sigma H$  \ \  is semisimple artinian;

(v) $R $ is left (right ) semi-hereditary iff  $R\#_\sigma H$ \ \ is
left (right ) semi-hereditary;

(vi)  $R $  is von Neumann regular iff $R\#_\sigma H $  \ \ is von Neumann regular.

\end {Corollary}

Since group algebra $kG$  is a cocommutative Hopf algebra, we have
that
     $$rgD(R) = rgD(R*G).$$
     Thus Corollary \ref {2.9} implies  in \cite
[Theorem 7.5.6] {MR}.

\vskip 2cm


\end {document}